\documentclass{paper}
\usepackage[utf8]{inputenc}
\usepackage{float}
\usepackage{amsmath}
\usepackage{amsthm}
\usepackage{amssymb}
\usepackage{graphicx}
\usepackage[numbers]{natbib}
\usepackage[unicode=true,
 bookmarks=false,
 breaklinks=false,pdfborder={0 0 1},backref=section,colorlinks=false]
 {hyperref}

\makeatletter

\theoremstyle{plain}
    \ifx\thechapter\undefined
      \newtheorem{prop}{\protect\propositionname}
    \else
      
    \fi
\theoremstyle{definition}
    \ifx\thechapter\undefined
      \newtheorem{defn}{\protect\definitionname}
    \else
      
    \fi

\usepackage{stmaryrd}

\makeatother

\providecommand{\definitionname}{Definition}
\providecommand{\propositionname}{Proposition}

\begin{document}
\title{Inertial control of a spinning flat disk}
\author{Luc Jaulin}
\institution{Lab-Sticc, ENSTA-Bretagne}

\maketitle
\textbf{Abstract}. This paper proposes a Lyapunov based approach to
control the rotation of a flat disk spinning in the space without
external forces. The motion of the disk is governed by the Euler's
rotation equation for spinning objects. The control is made through
the inertia matrix of the disk.

\section{Spinning disk}

Consider a disk spinning in the space without any gravity. The disk
is assumed to be flat, \emph{i.e.}, its inertia matrix \citep{wells67}
\begin{equation}
\mathbf{I}=\left(\begin{array}{ccc}
I_{1} & 0 & 0\\
0 & I_{2} & 0\\
0 & 0 & I_{3}
\end{array}\right)
\end{equation}
should satisfy 
\begin{itemize}
\item the positivity condition: $I_{1}>0,I_{2}>0,I_{3}>0$;
\item the flatness condition: $I_{1}=I_{2}+I_{3}$.
\end{itemize}
We assume that we can control the inertia matrix via some internal
forces in a strap down manner. It means that we have no inertial wheels.
Instead, we could slightly modify $\mathbf{I}$ by dilatation or compression
of some parts of the disk. We also consider that we are not able to
change the center of gravity of the disk. 

\begin{figure}[H]
\begin{centering}
\includegraphics[width=8cm]{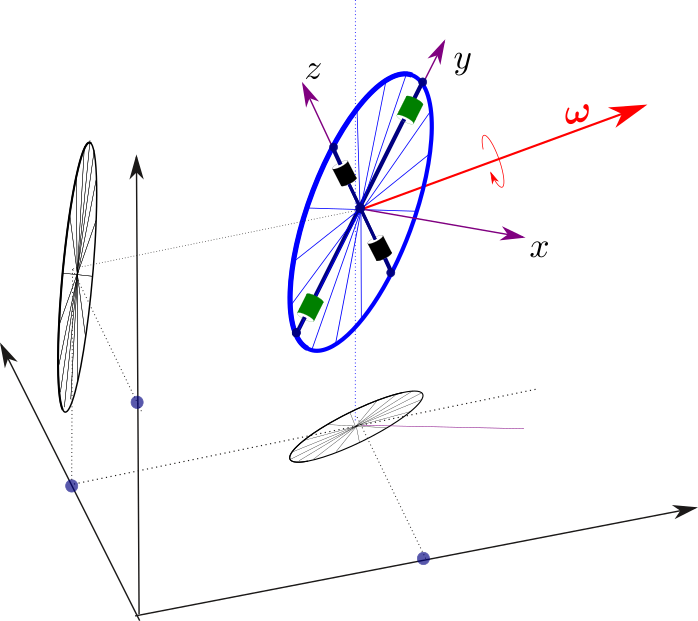}
\par\end{centering}
\caption{The disk spins and the two pairs of masses (black and green) can slightly
move in order the control the motion of the rotation vector $\boldsymbol{\omega}$ }
\label{fig:flatdisk}
\end{figure}

At rest, the inertial matrix is assumed to be

\begin{equation}
\bar{\mathbf{I}}=\left(\begin{array}{ccc}
\bar{I}_{1} & 0 & 0\\
0 & \bar{I}_{2} & 0\\
0 & 0 & \bar{I}_{3}
\end{array}\right)=\left(\begin{array}{ccc}
\frac{1}{2}mr^{2} & 0 & 0\\
0 & \frac{1}{4}mr^{2} & 0\\
0 & 0 & \frac{1}{4}mr^{2}
\end{array}\right)
\end{equation}
where $r$ is the radius of the disk and $m$ is its mass. The entries
of the inertia matrix can change as soon as the positivity and the
flatness conditions are satisfied. Equivalently, and for a better
visualization, we may assume that, we have two symmetric pairs of
masses which can move along the $y$ and $z$ axis (see Figure \ref{fig:flatdisk}).

The inertia of the disk is now
\begin{equation}
\mathbf{I}=\left(\begin{array}{ccc}
\frac{1}{2}mr^{2}+\text{\ensuremath{\frac{1}{2}}\ensuremath{\ensuremath{\ell_{2}^{2}}+\ensuremath{\frac{1}{2}\ell_{3}^{2}}}} & 0 & 0\\
0 & \frac{1}{4}mr^{2}+\text{\ensuremath{\text{\ensuremath{\frac{1}{2}}\ensuremath{\ensuremath{\ell_{3}^{2}}}}}} & 0\\
0 & 0 & \frac{1}{4}mr^{2}+\text{\ensuremath{\frac{1}{2}}\ensuremath{\text{\ensuremath{\ell_{2}^{2}}}}}
\end{array}\right)\label{eq:matrix:I}
\end{equation}
where $\ell_{2}$ corresponds to the distance of the $y$-masses (green)
to the center and $\ell_{3}$ is the distance of the $z$-masses (black)
to the center. Our objective is to propose a controller to modify
the rotation of the disk. The control objective is not so clear at
the moment since rotation cannot be controlled arbitrarily. Indeed,
the angular momentum remains constant and is thus non controllable.

This paper is organized as follows. Section \ref{sec:State-equation}
gives the state equations of the spinning disk. Section \ref{sec:Controller}
proposes different controllers to change behavior of the rotation
of the disk. Section \ref{sec:Conclusion} concludes the paper.

\section{State equation of the flat disk \label{sec:State-equation}}

The state equations are obtained from the conservation of angular
momentum. If $\mathbf{R}$ defines the orientation matrix of a rigid
body (for us the flat disk), then the angular momentum is

\begin{equation}
\mathcal{L}=\mathbf{R}\mathbf{I}\boldsymbol{\omega}_{r}
\end{equation}
where $\boldsymbol{\omega}_{r}$ is the rotation vector of the body
expressed in its own frame. 

\begin{figure}[H]
\begin{centering}
\includegraphics[width=10cm]{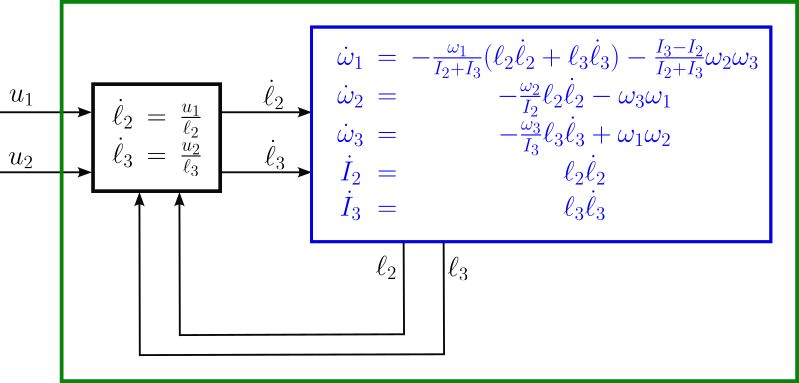}
\par\end{centering}
\caption{A first controller to simplify the state equations of the spinning
disk }
\label{fig:flatdiskloop1}
\end{figure}

Since it is constant, we have
\begin{equation}
\begin{array}{cl}
 & \dot{\mathbf{R}}\mathbf{I}\boldsymbol{\omega}_{r}+\mathbf{R}\dot{\mathbf{\cdot I}\cdot}\boldsymbol{\omega}_{r}+\mathbf{R}\cdot\mathbf{I}\cdot\dot{\boldsymbol{\omega}}_{r}=\mathbf{0}\\
\Leftrightarrow & \mathbf{R}^{\text{T}}\dot{\mathbf{R}}\mathbf{I}\boldsymbol{\omega}_{r}+\dot{\mathbf{I}}\cdot\boldsymbol{\omega}_{r}+\mathbf{I}\cdot\dot{\boldsymbol{\omega}}_{r}=\mathbf{0}
\end{array}
\end{equation}
We get the Euler's rotation equation
\begin{equation}
\dot{\boldsymbol{\omega}}_{r}=\mathbf{I}^{-1}\cdot\left(-\dot{\mathbf{I}}\cdot\boldsymbol{\omega}_{r}-\boldsymbol{\omega}_{r}\wedge\mathbf{I}\cdot\boldsymbol{\omega}_{r}\right)
\end{equation}
If we include the orientation, we have the motion equation
\begin{equation}
\left\{ \begin{array}{cclcl}
\mathbf{\dot{R}} &  & = &  & \mathbf{R}\cdot\textnormal{(}\boldsymbol{\omega}_{r}\wedge\textnormal{)}\\
\dot{\boldsymbol{\omega}}_{r} &  & = &  & \mathbf{I}^{-1}\cdot\left(-\dot{\mathbf{I}}\cdot\boldsymbol{\omega}_{r}-\boldsymbol{\omega}_{r}\wedge\mathbf{I}\cdot\boldsymbol{\omega}_{r}\right)
\end{array}\right.
\end{equation}
which will be used for the simulation. For the control we do not want
to control the orientation $\mathbf{R}$; we only want to control
$\boldsymbol{\omega}_{r}$. In what follows, the component of $\boldsymbol{\omega}_{r}$
will be denoted by $\omega_{1},\omega_{2},\omega_{3}$, \emph{i.e.},
we will write
\begin{equation}
\boldsymbol{\omega}_{r}=(\omega_{1},\omega_{2},\omega_{3})^{\text{T}}.
\end{equation}
Since $\mathbf{I}$ is diagonal, the time-varying Euler equation rewrites
into
\begin{equation}
\left\{ \begin{array}{ccc}
I_{1}\dot{\omega}_{1} & = & -\dot{I}_{1}\omega_{1}-(I_{3}-I_{2})\omega_{2}\omega_{3}\\
I_{2}\dot{\omega}_{2} & = & -\dot{I}_{2}\omega_{2}-(I_{1}-I_{3})\omega_{3}\omega_{1}\\
I_{3}\dot{\omega}_{3} & = & -\dot{I}_{3}\omega_{3}-(I_{2}-I_{1})\omega_{1}\omega_{2}
\end{array}\right.
\end{equation}
From (\ref{eq:matrix:I}), we have

\begin{equation}
\begin{array}{ccc}
I_{1} & = & \frac{1}{2}mr^{2}+\frac{1}{2}\ell_{2}^{2}+\frac{1}{2}\ell_{3}^{2}\\
I_{2} & = & \frac{1}{4}mr^{2}+\frac{1}{2}\ell_{2}^{2}\\
I_{3} & = & \frac{1}{4}mr^{2}+\frac{1}{2}\ell_{3}^{2}
\end{array}
\end{equation}
where we assumed that the sum of the control masses along the $y$axis
are equal to $1$. The same assumption is taken for the masses along
the $z$ axis. We have
\begin{equation}
\begin{array}{ccc}
\dot{I}_{1} & = & \ell_{2}\dot{\ell}_{2}+\ell_{3}\dot{\ell}_{3}\\
\dot{I}_{2} & = & \ell_{2}\dot{\ell}_{2}\\
\dot{I}_{3} & = & \ell_{3}\dot{\ell}_{3}
\end{array}
\end{equation}
The state equations of the system are thus
\begin{equation}
\left\{ \begin{array}{ccc}
\dot{\omega}_{1} & = & -\frac{\omega_{1}}{I_{2}+I_{3}}(\ell_{2}\dot{\ell}_{2}+\ell_{3}\dot{\ell}_{3})-\frac{I_{3}-I_{2}}{I_{2}+I_{3}}\omega_{2}\omega_{3}\\
\dot{\omega}_{2} & = & -\frac{\omega_{2}}{I_{2}}\ell_{2}\dot{\ell}_{2}-\omega_{3}\omega_{1}\\
\dot{\omega}_{3} & = & -\frac{\omega_{3}}{I_{3}}\ell_{3}\dot{\ell}_{3}+\omega_{1}\omega_{2}\\
\dot{I}_{2} & = & \ell_{2}\dot{\ell}_{2}\\
\dot{I}_{3} & = & \ell_{3}\dot{\ell}_{3}
\end{array}\right.
\end{equation}
A first simplification of the control problem can be obtained by adding
the following controller
\begin{equation}
\left\{ \begin{array}{ccc}
\dot{\ell}_{2} & = & \frac{u_{1}}{\ell_{2}}\\
\dot{\ell}_{3} & = & \frac{u_{2}}{\ell_{3}}
\end{array}\right.
\end{equation}
where $u_{1},u_{2}$ are some intermediate input (see Figure \ref{fig:flatdiskloop1}),
we get the state equations for the flat disk: 
\begin{equation}
\left\{ \begin{array}{ccc}
\dot{\omega}_{1} & = & -\frac{\omega_{1}}{I_{2}+I_{3}}(u_{1}+u_{2})-\frac{I_{3}-I_{2}}{I_{2}+I_{3}}\omega_{2}\omega_{3}\\
\dot{\omega}_{2} & = & -\frac{\omega_{2}}{I_{2}}u_{1}-\omega_{3}\omega_{1}\\
\dot{\omega}_{3} & = & -\frac{\omega_{3}}{I_{3}}u_{2}+\omega_{1}\omega_{2}\\
\dot{I}_{2} & = & u_{1}\\
\dot{I}_{3} & = & u_{2}
\end{array}\right.\label{eq:state:eq:I1I2I3}
\end{equation}

The state vector is $\mathbf{x}=(\omega_{1},\omega_{2},\omega_{3},I_{2},I_{3})$
and the input vector is $\mathbf{u}=(u_{1},u_{2})$.

\section{Controller\label{sec:Controller}}

This section proposes different controllers to change the spinning
behavior of the flat disk. The principle is to combine Lyapunov control
theory in the context of rotating object \citep{Leine21}.

\subsection{Principle}

To find a controller for our system (\ref{eq:state:eq:I1I2I3}), we
follow a Lyapunov control approach \citep{SontagBook}, \citep{KhalilHassa2002}.
For this, we choose a positive function $V(\mathbf{x})$, namely the
\emph{objective function}, such that $V(\mathbf{x})=0$ when the objective
is reached. 

Let us rewrite the state equation into an affine form \citep{Isidori95}:
\begin{equation}
\underset{=\dot{\mathbf{x}}}{\underbrace{\left(\begin{array}{c}
\dot{\omega}_{1}\\
\dot{\omega}_{2}\\
\dot{\omega}_{3}\\
\dot{I}_{2}\\
\dot{I}_{3}
\end{array}\right)}}=\underset{=\mathbf{f}(\mathbf{x})}{\underbrace{\left(\begin{array}{c}
-\frac{I_{3}-I_{2}}{I_{2}+I_{3}}\omega_{2}\omega_{3}\\
-\omega_{3}\omega_{1}\\
\omega_{1}\omega_{2}\\
0\\
0
\end{array}\right)}}+\underset{=\mathbf{g}_{1}(\mathbf{x})}{\underbrace{\left(\begin{array}{c}
-\frac{\omega_{1}}{I_{2}+I_{3}}\\
-\frac{\omega_{2}}{I_{2}}\\
0\\
1\\
0
\end{array}\right)}}\cdot u_{1}+\underset{=\mathbf{g}_{2}(\mathbf{x})}{\underbrace{\left(\begin{array}{c}
-\frac{\omega_{1}}{I_{2}+I_{3}}\\
0\\
-\frac{\omega_{3}}{I_{3}}\\
0\\
1
\end{array}\right)}}\cdot u_{2}.
\end{equation}

The time derivative of $V$ is
\begin{equation}
\dot{V}(\mathbf{x},\mathbf{u})=\mathcal{L}_{\mathbf{f}}V(\mathbf{x})+\mathcal{L}_{\mathbf{g}_{1}}V(\mathbf{x})\cdot u_{1}+\mathcal{L}_{\mathbf{g}_{2}}V(\mathbf{x})\cdot u_{2},
\end{equation}
where $\mathcal{L}$ is the Lie derivative operator. Since the solution
of 
\begin{equation}
\begin{array}{cc}
\min & a_{1}u_{1}+a_{2}u_{2}\\
\text{s.t.} & u_{1}^{2}+u_{2}^{2}=1
\end{array}
\end{equation}
where $\,a_{1}=\mathcal{L}_{\mathbf{g}_{1}}V(\mathbf{x}),\,a_{2}=\mathcal{L}_{\mathbf{g}_{2}}V(\mathbf{x})$,
is
\begin{equation}
\left(\begin{array}{c}
u_{1}\\
u_{2}
\end{array}\right)=-\frac{1}{\sqrt{a_{1}^{2}+a_{2}^{2}}}\cdot\left(\begin{array}{c}
a_{1}\\
a_{2}
\end{array}\right),
\end{equation}
we conclude that a Lyapunov like controller to reach our objective
is
\begin{equation}
\mathbf{u}=-\frac{1}{\sqrt{(\mathcal{L}_{\mathbf{g}_{1}}V(\mathbf{x}))^{2}+(\mathcal{L}_{\mathbf{g}_{2}}V(\mathbf{x}))^{2}}}\cdot\left(\begin{array}{c}
\mathcal{L}_{\mathbf{g}_{1}}V(\mathbf{x})\\
\mathcal{L}_{\mathbf{g}_{2}}V(\mathbf{x})
\end{array}\right).
\end{equation}

\subsection{Alignment controller\label{subsec:Alignment-controller}}

Assume that we want the disk spins around one principal axis of the
body. For instance its first axis, \emph{i.e.}, the $x$-axis. We
define the objective function
\begin{equation}
V(\mathbf{x})=\frac{1}{2}\left(\omega_{2}^{2}+\omega_{3}^{2}+(I_{2}-\bar{I}_{2})^{2}+(I_{3}-\bar{I}_{3})^{2}\right).
\end{equation}
The quantity $\omega_{2}^{2}+\omega_{3}^{2}$ corresponds to the \emph{alignment
error.} Note that when $V(\mathbf{x})=0$, we have $\omega_{2}=\omega$$_{3}=0$
(the alignment is performed) and $\mathbf{I}=\bar{\mathbf{I}}$ (the
inertia matrix is at its nominal position). The time derivative of
$V$ is
\begin{equation}
\begin{array}{ccl}
\dot{V}(\mathbf{x}) & = & \omega_{2}\dot{\omega}_{2}+\omega_{3}\dot{\omega}_{3}+(I_{2}-\bar{I}_{2})\dot{I}_{2}+(I_{3}-\bar{I}_{3})\dot{I}_{3}\\
                    & = & \omega_{2}\left(-\frac{\omega_{2}}{I_{2}}u_{1}-\omega_{3}\omega_{1}\right)+\omega_{3}\left(-\frac{\omega_{3}}{I_{3}}u_{2}+\omega_{1}\omega_{2}\right)
                          +  (I_{2}-\bar{I}_{2})u_{1}+(I_{3}-\bar{I}_{3})u_{2}\\
                    & = & -\frac{\omega_{2}^{2}}{I_{2}}u_{1}-\omega_{3}\omega_{1}\omega_{2}-\frac{\omega_{3}^{2}}{I_{3}}u_{2}+\omega_{1}\omega_{2}\omega_{3}+(I_{2}
                          -\bar{I}_{2})u_{1}+(I_{3}-\bar{I}_{3})u_{2}\\
                    & = & \underset{\mathcal{L}_{\mathbf{f}}V}{\underbrace{0\cdot\omega_{1}\omega_{2}\omega_{3}}}+\underset{\mathcal{L}_{\mathbf{g}_{1}}V}
                         {\underbrace{\left(-\frac{\omega_{2}^{2}}{I_{2}}+(I_{2}-\bar{I}_{2})\right)}}\cdot u_{1}+\underset{\mathcal{L}_{\mathbf{g}_{2}}V}
                         {\underbrace{\left(-\frac{\omega_{3}^{2}}{I_{3}}+(I_{3}-\bar{I}_{3})\right)}}\cdot u_{2}
\end{array}
\end{equation}

The same result would have easily been obtained using a symbolic calculus
module such as \texttt{Sympy} which is an open source computer algebra
system written in \texttt{Python}. The corresponding \texttt{Sympy}
script is the following

\begin{verbatim}
from sympy import*
from sympy.diffgeom import*
C=CoordSystem('C',Patch('P',Manifold('M',5)),{[w1,w2,w3,I2,I3]})     
w1,w2,w3,I2,I3 = C.coord\_functions()     
E = C.base\_vectors()
F=-((I3-I2)/(I2+I3))*w2*w3*E[0]+-w3*w1*E[1]+w1*w2*E[2]
G1=-(w1/(I2+I3))*E[0]-(w2/I2)*E[1]+E[3]
G2=-(w1/(I2+I3))*E[0]-(w3/I3)*E[2]+E[4]
V = 1/2*w2**2 + 1/2*w3**2 + 1/2*(I2-I20)**2 + 1/2*(I3-I30)**2
LfV=LieDerivative(F,V)
Lg1V=LieDerivative(G1,V)
Lg2V=LieDerivative(G2,V)
\end{verbatim}
\textbf{Stability analysis}. With our controller, we have $\dot{V}(\mathbf{x})=0$
if 
\begin{equation}
\left\{ \begin{array}{ccc}
-\frac{\omega_{2}^{2}}{I_{2}}+I_{2}-\bar{I}_{2} & = & 0\\
-\frac{\omega_{3}^{2}}{I_{3}}+I_{3}-\bar{I}_{3} & = & 0
\end{array}\right.
\end{equation}
We conclude that it it possible to converge to values $\bar{\omega}_{1},\bar{\omega}_{2},\bar{\omega}_{3}$
that do not correspond to the desired alignment. Equivalently, we
may have $\bar{\omega}_{2}\neq0,$ or $\bar{\omega}_{3}\neq0$. 

\begin{figure}[H]
\begin{centering}
\centering\includegraphics[width=12cm]{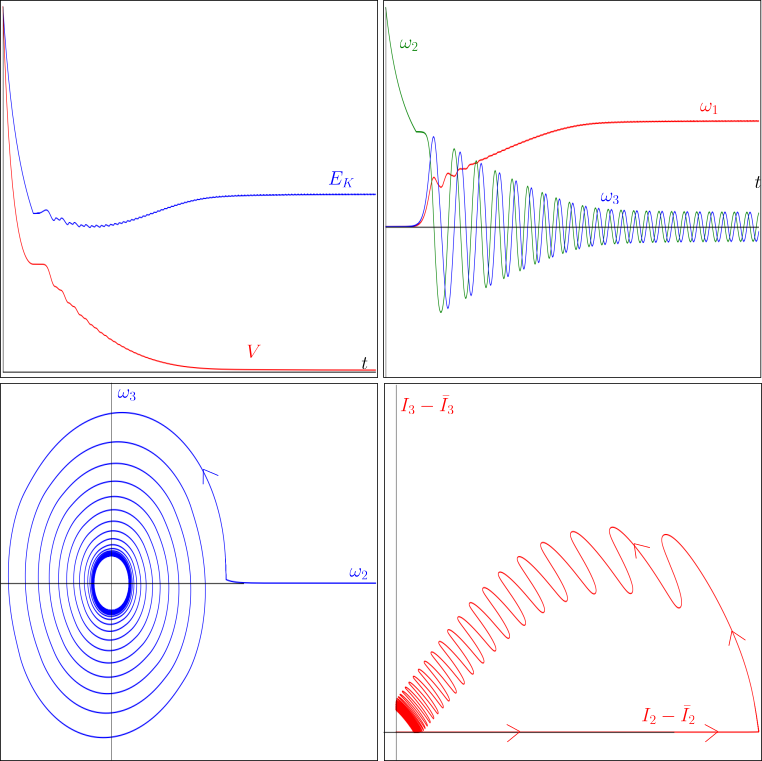}
\par\end{centering}
\caption{The alignment control tends to change radically the rotation vector
from one principal axis to another }
\label{fig:simu_align}
\end{figure}

Figure \ref{fig:simu_align} illustrates the behavior of the controller.
The initial rotation vector is $\boldsymbol{\omega}(0)=(10^{-5},10,0).$
The first subfigure (top-left) shows the kinetic energy $E_{K}$ of
the disk (blue) and the objective function $V(t)$ (red) with respect
to $t$ for $t\in[0,40]$. The second subfigure (top right) shows
the evolution of the components of $\boldsymbol{\omega}$ in the frame
box $[0,40]\times[-7,10].$ The third subfigure represents the evolution
of the vector $(\omega_{2},\omega_{3}).$ The frame box is $[-5,5]\times[-5,5]$.
The last subfigure represents the component $I_{2},I_{3}$ of the
inertia matrix in the frame box $[-0.2,2.5]\times[-0.2,2]$. Since
the controller is not able to cancel $(\omega_{2},\omega_{3})$, we
consider that the alignment objective is only partially reach.

The strategy followed by the controller is first to increase $I_{2}$
to make the second axis the intermediate axis (see the $4$th Subfigure).
As a consequence, due to instability of the second axis, also known
as the Dzhanibekov effect (even if first shown by Poinsot \citep{Poinsot}),
the component $\omega_{2}$ of $\boldsymbol{\omega}$ is able to decrease
without precession (see the $3$th Subfigure). Then in a second step,
the control changes its strategy to generate an oscillation between
$I_{2}$ and $I_{3}$. This yields a convergence of $\omega_{1}$,
the rotation along the first axis. We also observe that $\omega_{2}$
and $\omega_{3}$ do not converge to a static value.

\subsection{Passive controller}

We now want the disk losses energy with the intuition that it may
limit the precession. Another motivation is that it may help to find
an objective function which could help to guaranty the stability.
A passivity control approach can be used for this purpose \citep{fantoni:02}\citep{lozano:dissipative}.
We define the \emph{mechanical energy} as
\begin{equation}
\begin{array}{ccc}
V(\mathbf{x}) & = & \frac{1}{2}I_{1}\omega_{1}^{2}+\frac{1}{2}I_{2}\omega_{2}^{2}+\frac{1}{2}I_{3}\omega_{3}^{2}+\frac{1}{2}(I_{2}-\bar{I}_{2})^{2}+\frac{1}{2}(I_{3}-\bar{I}_{3})^{2}\end{array}.
\end{equation}
The quantity $\frac{1}{2}\boldsymbol{\omega}_{r}^{\text{T}}\mathbf{I}\boldsymbol{\omega}_{r}$
corresponds to kinetic energy and \emph{$\frac{1}{2}(I_{2}-\bar{I}_{2})^{2}+\frac{1}{2}(I_{3}-\bar{I}_{3})^{2}$
}is the artificial potential energy\emph{.} 

We have
\[
\begin{array}{ccl}
\dot{V}(\mathbf{x}) & = & I_{1}\omega_{1}\dot{\omega}_{1}+I_{2}\omega_{2}\dot{\omega}_{2}+I_{3}\omega_{3}\dot{\omega}_{3}+
                           \frac{1}{2}\omega_{1}^{2}\dot{I}_{1}+\frac{1}{2}\omega_{2}^{2}\dot{I}_{2}+\frac{1}{2}\omega_{3}^{2}\dot{I}_{3}\\
                    &    & +(I_{2}-\bar{I}_{2})\dot{I}_{2}+(I_{3}-\bar{I}_{3})\dot{I}_{3}\\
                    & =  & \left(I_{2}+I_{3}\right)\omega_{1}\left(-\frac{\omega_{1}}{I_{2}+I_{3}}(u_{1}+u_{2})
                           -\frac{I_{3}-I_{2}}{I_{2}+I_{3}}\omega_{2}\omega_{3}\right) \\
                    &    &    +I_{2}\omega_{2}\left(-\frac{\omega_{2}}{I_{2}}u_{1}- \omega_{3}\omega_{1}\right)
                           +  I_{3}\omega_{3}\left(-\frac{\omega_{3}}{I_{3}}u_{2}+\omega_{1}\omega_{2}\right)\\
                    &    & +\frac{1}{2}\omega_{1}^{2}(u_{1}+u_{2})+\frac{1}{2}\omega_{2}^{2}u_{1}+\frac{1}{2}\omega_{3}^{2}u_{2}
                           +(I_{2}-\bar{I}_{2})u_{1}+(I_{3}-\bar{I}_{3})u_{2}\\
 & = & -\frac{1}{2}\omega_{1}^{2}u_{1}-\frac{1}{2}\omega_{2}^{2}u_{1}-\frac{1}{2}\omega_{1}^{2}u_{2}-\frac{1}{2}\omega_{3}^{2}u_{2}+(I_{2}-\bar{I}_{2})u_{1}+(I_{3}-\bar{I}_{3})u_{2}\\
 & = & \underset{\mathcal{L}_{\mathbf{f}}V}{\underbrace{0}}+\underset{\mathcal{L}_{\mathbf{g}_{1}}V}{\underbrace{\left(-\frac{1}{2}\omega_{1}^{2}-\frac{1}{2}\omega_{2}^{2}+I_{2}-\bar{I}_{2}\right)}}\cdot u_{1}+\underset{\mathcal{L}_{\mathbf{g}_{2}}V}{\underbrace{\left(-\frac{1}{2}\omega_{1}^{2}-\frac{1}{2}\omega_{3}^{2}+I_{3}-\bar{I}_{3}\right)}}\cdot u_{2}
\end{array}
\]
The same result is obtained using \texttt{Sympy} taking

\begin{verbatim}
V = 1/2*(I2+I3)*w1**2+ 1/2*I2*w2**2 + 1/2*I3*w3**2 
    + 1/2*(I2-I20)**2 + 1/2*(I3-I30)**2
\end{verbatim}

\begin{figure}[H]
\begin{centering}
\centering\includegraphics[width=12cm]{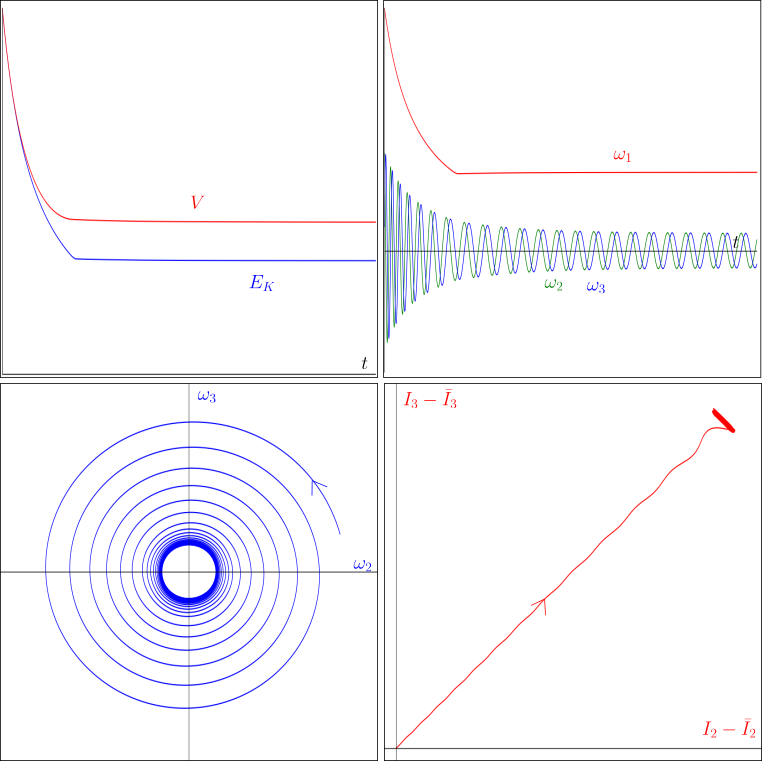}
\par\end{centering}
\caption{For the flat disk, a passivity approach is not sufficient to cancel
the precession}
\label{fig:simu_passiv}
\end{figure}

Figure \ref{fig:simu_passiv} illustrates the behavior of the controller
for $t\in[0,40]$. The initial rotation vector is $\boldsymbol{\omega}(0)=(10,4,1).$
The meaning of the plots are those given for Figure \ref{fig:simu_align}.
The frame for the subfigures are $[0,40]\times[0,1]$, $[0,40]\times[-5,10],$
$[-5,5]\times[-5,5]$ and $[-0.2,6]\times[-0.2,6]$, respectively.
We observe that the passivity approach for the flat disk is not sufficient
to cancel the precession. 

\subsection{Precession controller}

To cancel the precession, we define the objective function 
\begin{equation}
\begin{array}{ccc}
V(\mathbf{x}) & = & \frac{1}{2}\left((I_{3}-I_{2})\omega_{2}\omega_{3}\right)^{2}+\frac{1}{2}\left((I_{1}-I_{3})\omega_{3}\omega_{1}\right)^{2}+\frac{1}{2}\left((I_{2}-I_{1})\omega_{1}\omega_{2}\right)^{2}\\
 &  & +\frac{1}{2}(I_{2}-\bar{I}_{2})^{2}+\frac{1}{2}(I_{3}-\bar{I}_{3})^{2}.
\end{array}
\end{equation}
The quantity 
\begin{equation}
\left((I_{3}-I_{2})\omega_{2}\omega_{3}\right)^{2}+\left((I_{1}-I_{3})\omega_{3}\omega_{1}\right)^{2}+\left((I_{2}-I_{1})\omega_{1}\omega_{2}\right)^{2}=\|\boldsymbol{\omega}_{r}\wedge\mathbf{I}\cdot\boldsymbol{\omega}_{r}\|^{2}
\end{equation}
can be interpreted precession energy, \emph{i.e.}, a part of the kinetic
energy that creates the precession and that can be recovered. Indeed,
when $\boldsymbol{\omega}_{r}\wedge\mathbf{I}\cdot\boldsymbol{\omega}_{r}=\mathbf{0}$,
we see from the Euler equation that $\dot{\boldsymbol{\omega}}_{r}=\mathbf{0}$
when $\mathbf{u}=\mathbf{0},$which means that we have no more precession.

Using \texttt{Sympy} with

\begin{verbatim}
V = 1/2*((I3-I2)*w2*w3)**2 + 1/2*((I2)*w1*w3)**2
     + 1/2*((-I3)*w1*w2)**2
     + 1/2*(I2-I20)**2 + 1/2*(I3-I30)**2
\end{verbatim}

we get
\[
\begin{array}{ccc}
\mathcal{L}_{\mathbf{f}}V(\mathbf{x}) & = & \frac{\left(\left(\omega_{2}^{2}I_{3}^{2}+\omega_{3}^{2}I_{2}^{2}\right)\left(I_{2}-I_{3}\right)+\left(I_{2}+I_{3}\right)\left(\left(I_{2}-I_{3}\right)^{2}\omega_{2}^{2}-\left(I_{2}-I_{3}\right)^{2}\omega_{3}^{2}+\omega_{1}^{2}I_{2}^{2}-\omega_{1}^{2}I_{3}^{2}\right)\right)\omega_{1}\omega_{2}\omega_{3}}{I_{2}+I_{3}}\\
\mathcal{L}_{\mathbf{g}_{1}}V(\mathbf{x}) & = & \frac{1-\left(\left(I_{2}-I_{3}\right)^{2}\omega_{3}^{2}+\omega_{1}^{2}I_{3}^{2}\right)\omega_{2}^{2}}{I_{2}}-\frac{\left(\omega_{2}^{2}I_{3}^{2}+\omega_{3}^{2}I_{2}^{2}\right)\omega_{1}^{2}}{I_{2}+I_{3}}+\frac{\left(I_{2}-\bar{I}_{2}\right)+\left(I_{2}-I_{3}\right)\omega_{2}^{2}\omega_{3}^{2}+\omega_{1}^{2}\omega_{3}^{2}I_{2}}{I_{2}+I_{3}}\\
\mathcal{L}_{\mathbf{g}_{2}}V(\mathbf{x}) & = & \frac{1-\left(\left(I_{2}-I_{3}\right)^{2}\omega_{2}^{2}+\omega_{1}^{2}I_{2}^{2}\right)\omega_{3}^{2}}{I_{3}}-\frac{\left(\omega_{2}^{2}I_{3}^{2}+\omega_{3}^{2}I_{2}^{2}\right)\omega_{1}^{2}}{I_{2}+I_{3}}+\frac{\left(I_{3}-\bar{I}_{3}\right)+\left(I_{3}-I_{2}\right)\omega_{2}^{2}\omega_{3}^{2}+\omega_{1}^{2}\omega_{2}^{2}I_{3}}{I_{2}+I_{3}}
\end{array}
\]

\begin{figure}[H]
\begin{centering}
\centering\includegraphics[width=12cm]{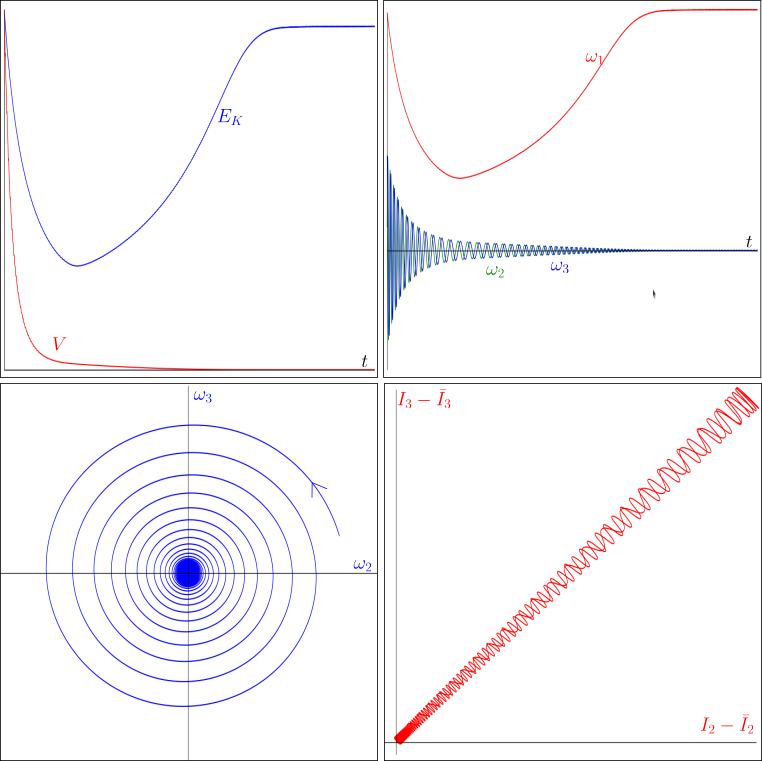}
\par\end{centering}
\caption{The controller cancel the precession efficiently}
\label{fig:simu_precession}
\end{figure}

Figure \ref{fig:simu_precession} illustrates the behavior of the
controller for $t\in[0,80]$. Again, the initial rotation vector is
$\boldsymbol{\omega}(0)=(10,4,1).$ The meaning of the plots are those
given for Figure \ref{fig:simu_align}. The frame for the subfigures
are $[0,80]\times[0,1]$, $[0,80]\times[-5,10],$ $[-5,5]\times[-5,5]$
and $[-0.2,6]\times[-0.2,6]$ respectively. We observe that the controller
easily cancel the precession of the flat disk.

\subsection{Combination}

As for Subsection \ref{subsec:Alignment-controller}, assume that
at time $t=0$, we spin around the second axis: 
\[
\boldsymbol{\omega}_{r}(0)=(10^{-5},10,0)
\]
and that we want to spin around the second axis. Then, it seems natural
to apply the alignment controller for a while and then to switch to
the precession controller to make $\boldsymbol{\omega}_{r}$ converge
to the nearest principal axis.

\begin{figure}[H]
\begin{centering}
\centering\includegraphics[width=12cm]{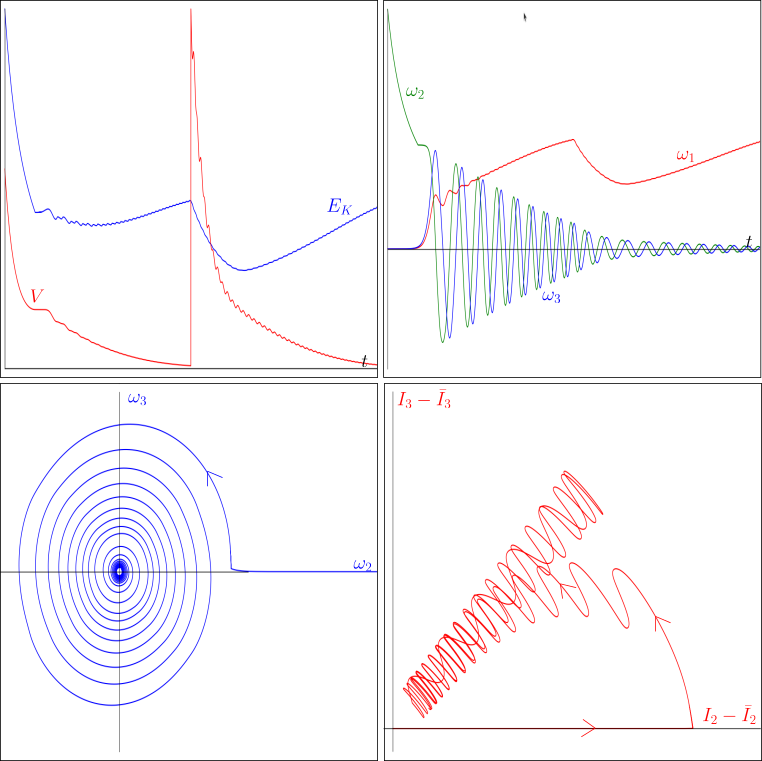}
\par\end{centering}
\caption{The combination of the two controllers solves the alignment problem
perfectly }
\label{fig:simu_combinaison}
\end{figure}

Figure \ref{fig:simu_combinaison} illustrates the behavior of the
controller. For $t\in[0,20]$, we apply the alignment control and
for $t\in[20,40]$, we apply the precession control. The objective
can be considered has reached since we were able to go from $\boldsymbol{\omega}_{r}=(10^{-5},10,0)$
to $\boldsymbol{\omega}_{r}=(5,0,0)$ without any external force.
This discontinuity observed for $V$ on the first subfigure is due
to the change of objective function at time $t=20.$ The frame for
the subfigures are $[0,40]\times[0,1]$, $[0,40]\times[-5,10],$ $[-5,5]\times[-5,5]$
and $[-0.2,4]\times[-0.2,3]$, respectively. 

A video illustrating the disk spinning and changing its rotation thanks
to the moving masses in the flat disk is given at: 
\begin{center}
\href{https://youtu.be/gRzlDYFuMts}{https://youtu.be/gRzlDYFuMts}
\par\end{center}

\section{Conclusion\label{sec:Conclusion}}

This paper has proposed several Lyapunov-like methods to control the
rotation of a spinning disk without any external forces. Only a partial
control of the rotation can be done since the angular momentum remains
constant. We made the assumption that the disk is flat which yields
a difficulty since masses are not allowed to move along the $x$ axis
of the disk. Even with this loss of controllability, we have shown
that some partial control was still feasible.

The first controller tries to cancel two components of the spin vector
to facilitate the alignment of the desired principal axis with the
angular momentum. Simulations show that the alignment is not perfectly
reached since some precession still exist. The second controller forces
the mechanical energy to decrease. The motivation of using a passive
control approach was to facilitate the proof of stability. Now, due
to the flatness assumption of the disk, the stability was not obtained
and we were not able to cancel the precession using an energy based
controller. The third controller focuses on the precession and is
shown to perfectly cancel it. This means that the spin vector will
perfectly be aligned with one of the principal axis of the body. Finally,
we have shown that the combination of the alignment controller with
the precession controller allows us to change the rotation axis from
one principal axis to another in an efficient way. Now, this combination
requires two objective functions whereas one should probably be sufficient.
A perspective of this work is to find a unique objective function
which would allow the controller to switch from one principal axes
to another. This objective function would probably give us a better
understanding of the motion of a flat disk. 

A possible application of the tools presented in this paper could
be for instance, the precession control of a flat rolling wheel \citep{Esnault_rapport2023}.
On the one hand, maintaining a precession can be helpful to change
rapidly the heading of the wheel. On the other hand, reducing the
precession increases the stability of the rolling wheel. 

\bigskip{}

The Python code associated to all examples can be found at
\begin{center}
\href{https://www.ensta-bretagne.fr/jaulin/flatdisk.html}{https://www.ensta-bretagne.fr/jaulin/flatdisk.html}
\par\end{center}

\bibliographystyle{plain}

\end{document}